\documentclass[12pt,leqno]{amsart}
\usepackage{amsmath,amsfonts,amssymb,amsthm}
\usepackage{vmargin}

\newtheorem{theorem}{Theorem}
\newtheorem{proposition}[theorem]{Proposition}
\newtheorem{lemma}[theorem]{Lemma}

\newtheorem{corollary}[theorem]{Corollary}

\theoremstyle{definition}
\newtheorem{remark}[theorem]{Remark}
\newtheorem{definition}[theorem]{Definition}
\newtheorem{example}[theorem]{Example}

\newcommand{\field}[1]{\mathbb{#1}}
\newcommand{\N}{\field{N}}                      
\newcommand{\F}{\field{F}}                      
\newcommand{\G}{\field{G}}                      
\newcommand{\R}{\field{R}}                      
\newcommand{\C}{\field{C}}                      
\newcommand{\Sph}{\field{S}}                      
\newcommand{\Heis}{\field{H}}                   

\newcommand{\cC}{\mathcal C}

\newcommand{\cL}{\mathcal L}
\newcommand{\cQ}{\mathcal Q}

\newcommand{\deriv}[1]{{\frac{\partial}{\partial #1}}}

\DeclareMathOperator{\spa}{span}

\numberwithin{theorem}{section}

\title[Helical CR structures and sub-Riemannian geodesics]{Helical
CR structures and sub-Riemannian geodesics}
\author{John P. D'Angelo \and Jeremy T. Tyson}
\date{\today}
\address{Mathematics Department, Univ.\ of Illinois, 1409 W. Green
  St., Urbana IL 61801.}
\email{tyson@math.uiuc.edu, jpda@math.uiuc.edu}
\thanks{JPDA supported by NSF grant DMS-0500765. JTT supported by NSF
  grant DMS-0555869.}

\begin{document}

\begin{abstract}
A helical CR structure is a decomposition of a real Euclidean
space into an even-dimensional horizontal subspace and its orthogonal
vertical complement, together with an almost complex
structure on the horizontal space and a marked vector in the vertical
space. We prove an equivalence between such structures and step two
Carnot groups equipped with a distinguished normal geodesic, and also
between such structures and smooth real curves whose derivatives have
constant Euclidean norm. As a consequence, we relate step two Carnot
groups equipped with sub-Riemannian geodesics with this family of
curves. The restriction to the unit circle of certain planar homogeneous
polynomial mappings gives an instructive class of examples. We
describe these examples in detail.
\end{abstract}

\maketitle

\section{Introduction}\label{sec1}

In this paper we identify a connection between CR geometry and
sub-Riemannian geometry. We introduce helical CR structures to relate
the structure of geodesics in step two Carnot groups with a certain
family of smooth real curves.

\begin{definition}\label{defn2}
A {\em helical CR structure ${\cC}$ of type $(2n,p)$} on a Euclidean
space ${\R}^d = \R^{2n+p}$ is an orthogonal decomposition
\begin{equation}\label{eqn35}
{\R}^d = {\R}^{2n} \oplus {\R}^p
\end{equation}
together with an invertible skew-symmetric operator $A$ on ${\R}^{2n}$
and a vector $w\in {\R}^p$. The {\em horizontal space} of $\cC$ is the
subspace $H = {\R}^{2n}$, while $w$ is the {\em vertical direction}.

A {\em marked helical CR structure} is a triple $({\cC},v,u_0)$,
where ${\cC}$ is a helical CR structure, $v\in{\R}^{2n}$, and
$u_0 = v_0 \oplus w_0$ is a vector in ${\R}^d = {\R}^{2n} \oplus
\R^p$. In this case, we also write $u=v \oplus w$, where $w$ is the
vector in ${\R}^p$ from the previous paragraph.
\end{definition}

We say that $\cC$ is (a) {\it horizontally trivial} if $n=0$, (b)
{\em vertically trivial} if $p=0$, and (c) {\em completely
nontrivial} if $n>0$ and $p>0$.

A {\it Carnot group} is a nilpotent stratified Lie group. Such
groups are naturally equipped with so-called {\it
Carnot-Carath\'eodory metrics}, which are singular metrics (from the
perspective of Riemannian geometry). Carnot groups figure
prominently as examples in the theory of analysis in metric measure
spaces. They also arise as local tangent space models for
sub-Riemannian manifolds. Geodesics in sub-Riemannian manifolds
yield solutions to various path planning problems in control theory
and related applications. See \cite{CDPT} and \cite{M} for more
information.

In section \ref{carnot} we make precise the equivalence between
completely nontrivial helical CR structures and the geometry of step
two Carnot groups and their (normal) sub-Riemannian geodesics. The
principal results of this section are Proposition \ref{A1} and Theorem
\ref{A2}. The proof of Theorem \ref{A2} is a straightforward 
application of the method of bicharacteristics to solve the
sub-Riemannian Hamiltonian equations. The explicit form of geodesics
in step two Carnot groups dates back to work of Gaveau \cite{Ga} and
Brockett \cite{B1}, \cite{B2}; see Montgomery's book \cite{M} for a
very readable summary. However, the precise geometric
characterizations of such groups and their geodesics given in
Proposition \ref{A1} and Theorem \ref{A2} enable us to link their
study with an ostensibly very different subject, namely, the study of
a family of smooth real curves motivated by homogeneous polynomial
mappings between balls and spheres. The first author has made
considerable use of such polynomial mappings in complex Euclidean
spaces. See \cite{D1}, \cite[Chapter 5]{D2}, \cite{D4}, \cite{DLP}.
The notion of a helical CR structure turns out to be equivalent to
real curves which are the analogues of such mappings. We thus obtain
an elementary differential geometric characterization of such
structures and in turn of step two Carnot groups and their geodesics.

More precisely, let us consider the class $\cQ_1$ consisting of
smooth curves taking values in a real Euclidean space all of whose
derivatives have constant norm, and the subset $\cQ_0$ consisting of
curves which have constant norm themselves (i.e., whose image lies
in a sphere centered at the origin). In Theorem \ref{thmX} we show
that each $\gamma \in \cQ_0$ determines a canonical decomposition
${\R}^d = {\R}^{2n} \oplus {\R}^p$ and an identification of its {\em
horizontal space} ${\R}^{2n}$ with ${\C}^n$. The appearance of a
complex structure on the horizontal space is an intriguing
phenomenon which is by no means evident from the original definition
of the family $\cQ_0$. It links the study of such curves with the
notion of helical CR structure. Indeed, each such curve is of the
form $\gamma(s) = \exp(As) v \oplus w$, where $A$ is skew-symmetric
and invertible. From a geometric perspective, such curves can be
viewed as generalized helices; from a different point of view, the
induced helical CR structure naturally determines a set of planes
(two-dimensional subspaces of the horizontal space) into which the
projection of the curve is a circle. Examples include classical
helices in $\R^3$ as well as the skew-line on the torus
$\Sph^1\times\Sph^1$. Theorem \ref{thmB} formalizes these
observations and summarizes the relationship between helical CR
structures and the classes $\cQ_0$ and $\cQ_1$.

By combining Theorems \ref{A2} and \ref{thmB} we obtain the
following differential geometric characterization of step two Carnot
groups of contact type and their geodesics. Theorem \ref{maintheorem}
is restated and proved later in the paper as Corollary \ref{corC}(b).

\begin{theorem}\label{maintheorem}
Each nonaffine curve $\mu \in \cQ_1$ contained in a hyperplane of
$\R^d$ determines a unique step two stratified Lie algebra of
contact type on a subspace of $\R^d$ together with the germ of a
normal geodesic $c$ for the induced Carnot-Carath\'eodory metric on
the associated Lie group. Conversely, each such Lie algebra and
geodesic determine a curve in $\cQ_1$. The horizontal projections of
$\mu$ and $c$ coincide.
\end{theorem}

The simplest example of a completely nontrivial helical CR structure
arises from the identification of ${\R}^3$ with ${\C}\times{\R}$;
the skew-symmetric $A$ is the standard almost complex structure matrix
\begin{equation}\label{J}
J = \begin{pmatrix} 0 & -1 \\  1 & 0 \end{pmatrix}
\end{equation}
and the vertical direction $w$ is $(0,0,1)$. The associated curve in
${\cQ}_1$ is a right cylindrical helix, while the associated step two
Carnot group is the Heisenberg group $\Heis^1$. The
sub-Riemannian geodesics in $\Heis^1$ are certain nonlinear helices
whose horizontal projections are circles in the horizontal space $\R^2
= \C$. See Example \ref{heis}.

From the point of view of CR geometry, the unit sphere in $\C^n$ is
a CR manifold of hypersurface type. Its sub-Riemannian geometry is
modelled on the Heisenberg group. Theorem \ref{maintheorem} refers
to step two Carnot groups of contact type, i.e., with one missing
direction. In Theorem \ref{thm4} we extend the correspondence in
Theorem \ref{maintheorem} to relate arbitrary step two Carnot groups
with $p$-dimensional center, equipped with $p$ distinguished normal
geodesics, with $p$-tuples of elements in ${\cQ}_1$ with common
horizontal space but linearly independent vertical directions. More
precisely, for $1 \le j \le p$ let $\mu_j$ denote elements of
${\cQ}_1$ taking values in ${\R}^d$ and whose horizontal spaces
coincide. When their vertical directions are linearly independent we
naturally obtain a step two Carnot group of type $(2n,p)$.
Conversely, given any step two Carnot group, the structure matrices
for its Lie algebra together with a basis for its vertical space
naturally determine a $p$-tuple of such curves.

In section \ref{sec4} we carefully discuss the specific $\cQ_0$ curves
which arose from consideration of group-invariant proper mappings. For
these curves the corresponding skew-symmetric matrices are bidiagonal.
We compute their spectra. In Remark \ref{DRZ} we mention a recent
paper concerning bidiagonalization of skew-symmetric matrices.

We mention several directions for future work. One may consider the
analogous situation where the starting point is certain homogeneous
CR mappings between odd-dimensional complex spheres, and try to
relate such mappings with harmonic maps between Carnot groups. See
\cite{CL} for a study of such maps with Heisenberg group target. In
the setting of this paper, groups of step two arise because the unit
sphere in $\C^n$ is strongly pseudoconvex. By considering CR maps
between CR manifolds with degenerate Levi forms, it might be
possible to consider nilpotent groups of higher step. See \cite{RS}
for the first decisive use of groups of higher step to prove
regularity results for differential operators arising in CR
geometry. The connection which we identify between CR and
sub-Riemannian geometry may also be of interest in the higher
codimension setting. See \cite{BG} for the construction of partial
normal forms for the defining equations of CR manifolds of finite
type in higher codimension. The finite type conditions in \cite{BG}
arise via iterated commutators of complex vector fields and their
conjugates. It is also natural to study smoothly varying helical CR
structures on the tangent spaces of a manifold.

\section{Step two Carnot groups and helical CR structures}\label{carnot}

In this section we establish the relationship between helical CR
structures and the sub-Riemannian geometry of step two Carnot groups
and their geodesics. In particular, we prove Proposition \ref{A1} and
Theorem \ref{A2}.

\subsection{Stratified Lie algebras and Lie groups}\label{cgroup}

Let ${\bf g}$ be a Lie algebra. We assume that ${\bf g}$ admits a
direct sum decomposition ${\bf g} = {\bf v}_1 \oplus \cdots \oplus
{\bf v}_s$. We set ${\bf v}_i=0$ for $i>s$. We say that $\bf g$ is
{\it stratified} if
$$ [{\bf v}_i,{\bf v}_j] = {\bf v}_{i+j} $$
for all $i,j\ge 1$. The integer $s$ is called the {\it step} of the
graded or stratified structure.

Let ${\G}$ be the Lie group corresponding to ${\bf g}$. We identify
${\bf g}$ with the left-invariant vector fields on $\G$. We say
that $\G$ is {\it graded} or {\it stratified} if its Lie algebra
enjoys the corresponding property. Stratified Lie groups are also
known as {\it Carnot groups}.

Let $\G$ be a step two Carnot group with Lie algebra ${\bf g} = {\bf
v}_1 \oplus {\bf v}_2$. We write $m=\dim {\bf v}_1$ and $p=\dim {\bf
v}_2$. Directions in ${\bf v}_1$ are called {\it horizontal}, while
directions in ${\bf v}_2$ are called {\it extra} (or {\it vertical})
directions. The relation $[{\bf v}_1,{\bf v}_1] = {\bf v}_2$ implies
that $1 \le p \le \binom{m}{2}$. We say that ${\bf g}$ or ${\G}$ is
of {\it contact type} if $p=1$. Some authors (see \cite{M}) reserve
this terminology for a stronger bracket-generating condition.

Let $X_1, \ldots,X_m$ be a basis for ${\bf v}_1$ and let $T_1,
\ldots,T_p$ be a basis for ${\bf v}_2$. For each $i,j$ there are
constants $c_{ij}^\alpha$ such that
$$ [X_i, X_j] = \sum_{\alpha =1}^p c_{ij}^\alpha T_\alpha. $$
The skew-symmetric matrices $C^\alpha = (c_{ij}^\alpha)$ are called
the {\it structure matrices} of $\G$.

The exponential map is a global diffeomorphism from ${\bf g}$
to $\G$. We identify ${\bf g}$ with ${\R}^{m+p}$. In the
corresponding {\it exponential coordinates}
$(x_1,\ldots,x_{m};t_1,\ldots,t_p) \in {\R}^{m} \oplus {\R}^p$
we write
$$  X_i = \deriv{x_i} + \frac12 \sum_{j=1}^{m}
\sum_{\alpha=1}^p c_{ij}^\alpha x_j \deriv{t_\alpha} $$ for
$i=1,\ldots,m$, and $T_\alpha = \deriv{t_\alpha}$ for
$\alpha=1,\ldots,p$.

We consider briefly two extreme cases. When $p=1$ and $m=2n$ we
obtain the Heisenberg group ${\Heis}^n$; this is the unique Carnot
group of dimension $2n+1$ with one vertical direction. When $p =
\binom{m}{2}$ we obtain the so-called {\it free} two step nilpotent
Lie group ${\N}_{m,2}$. In this case, the number $p$ of vertical
directions is maximal for a given number of horizontal directions~$m$.

\begin{proposition}\label{A1}
A completely nontrivial helical CR structure on $\R^d$ uniquely
determines and is uniquely determined by a step two stratified Lie
algebra of contact type on a subspace of $\R^d$.
\end{proposition}

\begin{proof}
Let ${\cC}$ be a completely nontrivial helical CR structure on
${\R}^d$ with $2n,p,A=(a_{ij}),w$ as in Definition \ref{defn2}. Define
vector fields $X_1,\ldots,X_{2n}, T_w$ on ${\R}^{2n+p}$ by
$$ T_w = \sum_{\alpha = 1}^p w_\alpha \deriv{t_\alpha} $$
and
\begin{equation}\label{eqn36}
X_i = \deriv{x_i} + \frac12 \sum_{j=1}^{2n} a_{ij} x_j T_w.
\end{equation}
Thus $T_w$ is the constant vector field in ${\R}^{2n+p}$ pointing in
the direction of $w=(w_1,\ldots,w_p)$. The commutation relations
$[X_i,X_j] = a_{ij} T_w$ for $i,j=1,\ldots,2n$ and $[X_i, T_w] = 0$
for $i=1,\ldots,2n$ follow easily from \eqref{eqn36}. They
equip the vector subspace ${\R}^{2n}\oplus {\rm span}(T_w)$
with the structure of a step two stratified Lie algebra of contact
type.

Conversely, assume that ${\bf g} = {\bf v}_1 \oplus {\bf v}_2$ is a
subspace of ${\R}^d$ equipped with the structure of a step two
stratified Lie algebra of contact type. Let $B=(b_{ij})$ be the
structure matrix for ${\bf g}$, defined by the commutation relations
$[X_i, X_j] = b_{ij} T$ where $X_1,\ldots,X_{m}$ is a basis for the
first layer ${\bf v}_1$ and $T$ generates the second layer ${\bf
  v}_2$. Let $w$ span the one-dimensional space $\exp(T)$. The
skew-symmetric matrix $B$ need not be invertible, and hence cannot
necessarily be used in Definition \ref{defn2}. We obtain an invertible
skew-symmetric $A$ by restricting $B$ to the orthocomplement of its
null space. This space is even dimensional because nonzero eigenvalues
come in conjugate pairs; thus $A$ is a $(2n)\times(2n)$ matrix for
some $n$. This information now determines an orthogonal decomposition
${\R}^d = {\R}^{2n} \oplus {\R}^p$, an appropriate invertible
skew-symmetric $A$, and $w \in {\R}^p$, i.e., a helical CR structure
${\cC}$. Since $w \ne 0$, ${\cC}$ is vertically nontrivial. Since
${\bf g}$ is of step two, ${\cC}$ is horizontally nontrivial.
\end{proof}

\subsection{Carnot-Carath\'eodory metric and sub-Riemannian
  geodesics}\label{srgeodesics}

Let $\G$ be a step two Carnot group. We introduce a
sub-Riemannian metric $g_0$ as an inner product on the horizontal
distribution ${\bf v}_1$ by declaring the vector fields
$X_1,\ldots,X_{m}$ to be orthonormal. Such a metric defines a global
distance function on $\G$ by the formula
\begin{equation}\label{eqn37}
d_0(P,Q) = \inf {\rm length}_{g_0}(\gamma),
\end{equation}
where the infimum is taken over all smooth curves
$\gamma:[a,b]\to{\G}$ joining $P$ to $Q$ whose tangent vector
$\gamma'(s)$ lies in the subspace of $T_{\gamma(s)}{\G}$ spanned by
the vector fields $X_1,\ldots,X_{m}$. The bracket-generating
condition $[{\bf v}_1,{\bf v}_1] = {\bf v}_2$ guarantees that $d_0$
is a metric. See \cite{B} or \cite{M}. A curve $c$ that realizes the
infimum in \eqref{eqn37} is called a {\it sub-Riemannian geodesic}
joining $P$ to $Q$. Existence and regularity theory
for sub-Riemannian geodesics in Carnot groups is the focus of the
papers \cite{S1} and \cite{S2}.

A curve $c$ in $\G$ is called a {\it horizontal lift} of a curve
$\gamma:[a,b]\to\exp{\bf v}_1$ if $\pi\circ c = \gamma$, where
$\pi:{\G} \to \exp {\bf v}_1$ denotes the projection
$\pi(x_1,\ldots,x_m;t_1,\ldots,t_p) = (x_1,\ldots,x_m)$. By the
general theory of fiber bundles, for each curve $\gamma$ as above
and each $P\in \pi^{-1}(\gamma(a))$, there exists a unique
horizontal lift $c$ of $\gamma$ with $c(a) = P$.

A fundamental feature of sub-Riemannian geometry which distinguishes
it from its Riemannian counterpart is the potential existence of
(strictly) abnormal geodesics. A sub-Riemannian geodesic $c$ is
called {\it normal} if it is the projection to $\G$ of a
bicharacteristic curve in $T^*\G$, i.e., a solution to the
sub-Riemannian Hamiltonian system.  The book \cite{M} carefully
summarizes the delicate issues related to the regularity of
sub-Riemannian geodesics and especially to the existence of abnormal
minimizers. In the proof of the following theorem we make use of the
Hamiltonian formalism, which explains the restriction to normal
geodesics. Normality of geodesics in step two Carnot groups is studied
in the paper \cite{gk}.

\begin{theorem}\label{A2}
A completely nontrivial marked helical CR structure on $\R^d$
uniquely determines and is uniquely determined by a pair $(S,c)$,
where $S$ is a step two stratified Lie algebra of contact type on a
subspace of $\R^d$ and $c$ is the germ of a normal geodesic for the
sub-Riemannian (Carnot-Carath\'eodory) metric on the associated
Carnot group.
\end{theorem}

\begin{proof}
Let $({\cC},v,u_0)$ be a completely nontrivial marked helical CR
structure on ${\R}^d$. As in the proof of part (a) we introduce a
step two stratified Lie algebra of contact type on a subspace of
${\R}^d$. Let $X_i$ be the vector fields defined in \eqref{eqn36}.
Let $g_0$ be the sub-Riemannian metric on $\G$ for which the vector
fields $X_i$ are an orthonormal basis for the horizontal space.
Normal geodesics in $\G$ can be identified by the method of
bicharacteristics, projecting to $\G$ the solutions of Hamilton's
equations for the sub-Riemannian Hamiltonian. We use coordinates
$(x_1,\ldots,x_{2n},t)$ for points in $\G$, where $t$ denotes the
coefficient of $w$ in the vertical direction. The Hamiltonian is
$$
H(x,t;\xi,\tau) = \frac12  \sum_i \zeta_i^2,
$$
where $\zeta_i = \xi_i + \frac12 \sum_{j=1}^{2n} a_{ij} x_j \tau$
is dual to $X_i$. The bicharacteristics $(x(s),t(s);\xi(s),\tau(s))$
are solutions to Hamilton's equations $\dot{x} =
\nabla_\xi H$, $\dot{t} = \partial H/\partial\tau$, $\dot\xi = -
\nabla_x H$, $\dot\tau = -\partial H/\partial t$. Explicitly,
$$
\dot{x}_i = \zeta_i, \quad \dot{t} = \frac12 \sum_{j} a_{ij} x_j
\zeta_i, \quad \dot\xi_j = - \frac12 \sum_{i} a_{ij} \zeta_i
\tau, \quad \dot\tau = 0
$$
with initial conditions $x(0)=x_0$, $t(0)=t_0$, $\xi(0)=\xi_0$ and
$\tau(0) = \tau_0$. Thus $\tau(s) = \tau_0$ is constant in $s$.
Hamilton's equations can be written in the compact form
$$
\dot{x} = \zeta, \qquad \dot{t} = \frac12 \zeta^T  A  x, \quad
\dot\xi = \frac{\tau_0}2 A  \zeta, \quad \zeta = \xi + \frac{\tau_0}2
A  x.
$$
Thus $\dot\zeta = \tau_0 A \zeta$ and the curve $\zeta(s) = \exp(s
\tau_0 A) \zeta(0)$ is in ${\cQ}_0$. By an affine
reparameterization we may assume that either $\tau_0 = 0$ or
$\tau_0=1$.

We next obtain the following constant coefficient coupled ODE system
in ${\R}^{2n}$ for the first-layer position and momentum
variables $x,\xi$:
$$
\dot{x} = \frac{\tau_0}2 A  x + \xi, \quad \dot\xi = \frac{\tau_0^2}4
A^2  x + \frac{\tau_0}2 A  \xi,
$$
which we easily integrate to find
$$
x(s) = \frac12 ( \exp(sA) + I ) x_0 + ( \exp(sA) - I ) A^{-1}
\xi_0,
$$
$$
\xi(s) = \frac14 ( \exp(sA) - I) A x_0 + \frac12 ( \exp(sA) +
I ) \xi_0
$$
if $\tau_0 = 1$, or
$$
x(s) = x_0 + s \xi_0, \qquad \xi(s) = \xi_0
$$
if $\tau_0 = 0$. Note that in the case $\tau_0 = 1$ we may write
\begin{equation}\label{eqn38}
x(s) = \left( \frac12 x_0 - A^{-1} \xi_0 \right) + \exp(sA)
\left( \frac12 x_0 + A^{-1} \xi_0 \right).
\end{equation}
In either case we immediately see that $x(s)$ is the projection to
${\R}^{2n}$ of a curve in ${\cQ}_1$ obtained as the integral
of the curve $\gamma(s) = \exp (sA) v \oplus w$ for $v = \xi_0 +
\frac12 x_0$.

Conversely, assume that the pair $(S,c)$ consists of a subspace $S$
of ${\R}^d$ equipped with a step two stratified Lie algebra of
contact type and a normal geodesic $c:(-\epsilon,\epsilon) \to
(\exp(S),g_0)$, where $g_0$ denotes the Carnot-Carath\'eodory metric
on $\exp(S)$. As above, we write $c(s)=(x(s),t(s))$ and denote by
$(\xi(s),\tau(s))$ the corresponding momenta variables. As in the
proof of part (a), we associate to this data a completely nontrivial
helical CR structure ${\cC}$ on ${\R}^d$. To identify the marking
$({\cC},v,u_0)$, we set $u_0 = c(0)$. The discussion in the previous
paragraph shows that the projection of $c$ in the horizontal space
$\exp{\bf v}_1$ takes the form \eqref{eqn38} for suitable $x_0$ and
$\xi_0$. The proof is complete upon setting $v = \frac12 A x_0 +
\xi_0$.
\end{proof}

\begin{example}\label{heis}
The simplest step two Carnot group is the Heisenberg group
${\Heis}^1$, with Lie algebra ${\bf h}^1 = {\bf v}_1 \oplus {\bf
v}_2 = \spa\{X,Y\} \oplus \spa\{T\}$ and a single nontrivial bracket
relation $[X,Y]=T$. Sub-Riemannian geodesics in ${\Heis}^1$ take the
form
$$
s\mapsto (a + b e^{-is} , c + |b|^2 s - {\rm Im} ({ \overline a} b
e^{-is} ) ),
$$
see e.g.\ \cite{B} or \cite{CDPT}. The case $b=-a$, $c=0$
corresponds to geodesics issuing from the identity (origin) in
$\Heis^1 = \R^3$:
$$
s \mapsto (a - a e^{-is}, |a|^2 (s - \sin s)).
$$
The associated ${\cQ}_1$ curves from Theorem \ref{maintheorem} are
right cylindrical helices of the form
$$
s \mapsto (a - a e^{-is}, c + s)
$$
with $a \in {\C}$ and $c \in {\R}$. Note that all geodesics in
$\Heis^1$ (or more generally, in any step two Carnot group
satisfying the strong bracket-generating assumption) are normal.
\end{example}

\section{A characterization of
  helical CR structures}\label{diffgeom}

In this section, we show that the notion of helical CR structure can
be recovered from an elementary concept in real differential
geometry. We will study smooth real Euclidean curves all of whose
derivatives have constant norm. Such curves can be thought of as
``generalized helices''. We begin by developing some basic
properties of such curves. For instance, we show that the
G-curvatures (coupling parameters in the generalized Frenet
formulas) of such a curve are all constant in time. The main point
of the section is that every such curve induces a natural
decomposition of its target space into horizontal and vertical
directions thereby leading to a helical CR structure. See Theorem
\ref{thmX}.

\subsection{Curves with derivatives of constant norm}

A {\em smooth real curve} in ${\R}^d$ is a smooth map $\gamma$ from
an interval in ${\R}$ to ${\R}^d$. We say that $\gamma$ {\it lies
in} a set $S$ if its image is a subset of $S$. For $j\ge 0$, the
$j$th order derivative of $\gamma$ will be written $D^j\gamma$. We
let $\langle v,w \rangle$ denote the Euclidean inner product of
vectors $v$ and $w$ in ${\R}^d$, and $||v||^2$ denote the
corresponding squared norm. We use the same notation without regard
to the dimension $d$. Let $O(d)$ denote the orthogonal group acting
on $\R^d$.

\begin{definition}\label{defn1}
Let $k\ge 0$. We denote by ${\cQ}_k$ the collection of smooth real
curves $\gamma$ such that $D^j \gamma$ has constant Euclidean norm
whenever $j \ge k$.
\end{definition}

Thus $\gamma \in {\cQ}_0$ if for $j \ge 0$ there are constants $c_j$
such that
\begin{equation}\label{eqn2}
||D^j \gamma(s)||^2 = c_j
\end{equation}
for all $s$. Using $j=0$ in \eqref{eqn2} we see that
$\gamma$ lies in a sphere centered at the origin. On the
other hand, elements of ${\cQ}_1$ need not have bounded image.

For any skew-symmetric $A$ and $v\in\R^d$, the curve $\gamma$ defined
by $\gamma(s) = \exp(As) v$ is in ${\cQ}_0$. Indeed,
since $A$ commutes with $\exp(sA)$ and $||\exp(As)v||^2 = ||v||^2$ for
every $s$, we obtain
$$
||D^k \gamma(s)||^2 = ||\exp(As) A^k v||^2  = ||A^k v||^2.
$$
Thus $||D^k\gamma(\cdot)||$ is constant. In particular, if $\gamma$
satisfies the ODE $D\gamma = A\gamma$ for some skew-symmetric $A$,
then $\gamma\in\cQ_0$. More generally, if $A$ is skew-symmetric on 
${\R}^d$, $B\in O(d)$, and $v \in {\R}^d$, then the curve $\gamma$
defined by $\gamma(s) = B \exp(As) v $ is in ${\cQ}_0$. In fact,
$\gamma(s) = \exp(Cs) \gamma(0)$ where $C=BAB^{-1}$ is
skew-symmetric. Theorem \ref{thmX} provides a converse to these
assertions.

\begin{remark}
A smooth curve in $O(d)$ need not be a member of $\cQ_0$. For example,
put $\gamma(s) = A(s)v$, where $v$ is a unit vector in $\R^2$ and
$$
A(s) = \begin{pmatrix} \cos f(s) & -\sin f(s) \\ \sin f(s) & \cos
  f(s) \end{pmatrix}.
$$
Then $\gamma \in {\cQ}_0$ if and only if $f$ is affine. In general,
membership in $\cQ_0$ is preserved only under affine
reparameterization. See Remark \ref{rem1} for a discussion of the
issues concerning reparameterization.
\end{remark}

We discuss some operations which preserve the class ${\cQ}_0$. The
ideas here are related to the work in  
\cite{D2} concerning polynomial CR mappings between spheres of
different dimensions. Parts (i) and (iii) of Lemma \ref{lem1} show
that $\cQ_0$ is closed under orthogonal postcomposition and tensor
products. Part (ii) exhibits homotopies
joining arbitrary elements of $\cQ_0$ provided the target dimension
is increased appropriately.

\begin{lemma}\label{lem1}
Let $\gamma:{\R} \to {\R}^n$ and $\mu : {\R} \to {\R}^k$ be elements
of ${\cQ}_0$, let $T\in O(n)$, and let $\theta\in\R$. Then
\begin{enumerate}
\item[(i)] $T\circ \gamma \in {\cQ}_0$,
\item[(ii)] the $\theta$-juxtaposition $J_\theta(\gamma, \mu) =
(\cos\theta) \gamma \oplus (\sin\theta)\mu$ is in ${\cQ}_0$,
\item[(iii)] the tensor product $\gamma \otimes \mu$ is in ${\cQ}_0$.
\end{enumerate}
\end{lemma}

\begin{proof}
The proofs of (i) and (ii) are simple computations:
\begin{equation}\label{eqn7}
||D^k (T \circ \gamma)(s)||^2 = ||T \circ D^k \gamma(s)||^2 = ||D^k
\gamma(s)||^2 = c_k
\end{equation}
and
\begin{equation*}\begin{split}
||D^k &J_\theta(\gamma, \mu)(s)||^2
= ||J_\theta(D^k \gamma, D^k \mu)(s)||^2
= || (\cos\theta)D^k \gamma (s) \oplus (\sin\theta) D^k \mu (s)||^2
\\
&= || (\cos\theta)D^k \gamma (s)||^2 + || (\sin\theta) D^k \mu (s)||^2
= \cos^2\theta ||D^k \gamma(s)||^2 + \sin^2\theta ||D^k \mu(s)||^2.
\end{split}\end{equation*}
The latter expression is independent of $s$ since
$\gamma,\mu\in\cQ_0$.

Finally, we prove (iii). Let $f$ and $g$ be arbitrary curves in
Euclidean spaces of possibly different dimensions. By definition
\begin{equation}\label{eqn8}
\langle f \otimes g, u \otimes v \rangle = \langle f, u \rangle \
\langle g, v \rangle
\end{equation}
and in particular
\begin{equation}\label{eqn9}
||f \otimes g||^2 = ||f||^2 ||g||^2.
\end{equation}
It follows from \eqref{eqn9} that
$\frac{d}{dt} ||f \otimes g||^2 = \frac{d}{dt} ||f||^2 ||g||^2 =
2 \langle f, f'\rangle ||g||^2 + ||f||^2 2 \langle g, g'\rangle$.
Note that if $||h||^2$ is constant, then
\begin{equation}\label{eqn11}
0 = \frac{d}{dt} ||h||^2 = 2 \langle h, h' \rangle.
\end{equation}
Applying \eqref{eqn11} for both $h=f$ and $h=g$, and using
\eqref{eqn8} and \eqref{eqn9} yields
\begin{equation}\label{eqn111}\begin{split}
||\frac{d}{dt} &(f \otimes g)||^2 =|| f' \otimes g + f \otimes
g'||^2 \\ &= || f' \otimes g||^2 + 2 \langle f', f \rangle \langle g,
g' \rangle + ||f \otimes g'||^2 = ||f'||^2 ||g||^2 + ||f||^2 ||g'||^2.
\end{split}\end{equation}
Thus $f\otimes g$ and its first derivative have constant norm
whenever $f,g\in\cQ_0$. An easy induction on the number of
derivatives now yields the conclusion in (iii).
\end{proof}

The classes ${\cQ}_0$ or ${\cQ}_1$ have connections with generalized
Frenet frames and higher curvatures. Gluck \cite{Gl1}, \cite{Gl2}
considered smooth curves in ${\R}^n$ whose first $n$ derivatives are
linearly independent at each point. Applying Gram-Schmidt yields an
orthonormal frame along the curve; by differentiating this frame
field, one obtains higher curvatures (so-called G-curvatures), which
depend only upon the inner products $e_{kl}(s) = \langle D^k
\gamma(s), D^l \gamma (s) \rangle$ for $k,l \le n$. The following
lemma shows that all such inner products $e_{kl}(s)$ with $k,l\ge 1$
are independent of $s$ for a curve $\gamma \in \cQ_1$. We will use
Corollary \ref{lem3} in the proof of our characterization theorem
(Theorem \ref{thmX}). It also implies that all G-curvatures are
constant in the time parameter for such curves.

\begin{lemma}\label{lem2}
Let $\gamma \in \cQ_0$. Then all inner products $e_{kl}(s)$ are
independent of $s$ and $e_{kl} = 0$ when $k-l$ is odd. When $\gamma
\in {\cQ}_1$, the same conclusions hold for all $k,l \ge 1$.
\end{lemma}

\begin{corollary}\label{lem3}
Let $\gamma \in {\cQ}_0$. Let $W_e(s)$ denote the span of the
derivatives of $\gamma$ of even order at $s$ (including $\gamma(s)$
itself), and let $W_o(s)$ denote the span of the derivatives of
$\gamma$ of odd order at $s$. Then, for each $s$, the spaces $W_e(s)$
and $W_o(s)$ are orthogonal.
\end{corollary}

\begin{proof}[Proof of Lemma \ref{lem2}]
Assume \eqref{eqn2} holds. For each $k\ge 0$ we have
\begin{equation}\label{eqn12}
0 = D \left( ||D^k \gamma(s)||^2\right) = 2 \langle D^k \gamma(s),
D^{k+1} \gamma(s) \rangle.
\end{equation}
By \eqref{eqn12} the conclusion holds when $l=k+1$. Differentiating again shows
that
\begin{equation}\label{eqn13}
0 = ||D^{k+1} \gamma(s)||^2 +  \langle D^k \gamma(s), D^{k+2}
\gamma(s) \rangle.
\end{equation}
The squared norm term in \eqref{eqn13} is a constant, and hence so is the
inner product term. Therefore $e_{kl}$ is a constant when $l=k+2$.
By differentiating \eqref{eqn13} and using induction we obtain the first
statement for ${\cQ}_0$. By differentiating \eqref{eqn12} an even number
of times and using induction we obtain the second statement for
${\cQ}_0$. The statement for ${\cQ}_1$ follows because $\gamma
\in {\cQ}_1$ implies $D \gamma \in {\cQ}_0$. \end{proof}

\subsection{A characterization of ${\cQ}_0$ and ${\cQ}_1$}\label{sec3}

In Lemma \ref{lem4} we show that each element $\gamma \in {\cQ}_0$
satisfies a constant coefficient ODE, whose characteristic equation
has only purely imaginary solutions. From this we deduce Theorem
\ref{thmX}, which demonstrates how $\gamma$ induces a helical CR
structure in the target space.

\begin{lemma}\label{lem4}
Let $\gamma \in {\cQ}_0$. Then there is a positive integer $k$ and
constant scalars $C_j$ for $0 \le j < k$ such that
\begin{equation}\label{eqn14}
D^{2k} \gamma = \sum_{j=0}^{k-1} C_j D^{2j} \gamma.
\end{equation}
\end{lemma}

\begin{proof}
Since the target space is finite dimensional, for each $s$ there is a
smallest $k$ so that a non-trivial dependence among the even
derivatives of $\gamma$ holds in a neighborhood of $s$:
\begin{equation}\label{eqn15}
D^{2k} \gamma(s) = \sum_{j=0}^{k-1} C_j(s)  D^{2j} \gamma(s).
\end{equation}
Observe that the coefficient functions $C_j$ in \eqref{eqn15} are smooth.
Differentiating \eqref{eqn15} gives
\begin{equation}\label{eqn16}
D^{2k+1} \gamma(s) = \sum_{j=0}^{k-1} C_j(s)  D^{2j+1} \gamma(s) +
\sum_{j=0}^{k-1} C_j'(s) D^{2j} \gamma (s).
\end{equation}
The spaces $W_e(s)$ and $W_o(s)$ are orthogonal by Corollary \ref{lem3}.
Orthogonality of the appropriate terms in \eqref{eqn16} forces
\begin{equation}\label{eqn17}
0 = \sum_{j=0}^{k-1} C_j'(s) D^{2j}\gamma(s).
\end{equation}
Since \eqref{eqn17} gives a linear dependence among the even derivatives for
a smaller value of $k$, and $k$ was chosen minimally, relation \eqref{eqn17}
must be trivial. Thus $C_j'(s)$ vanishes for each $j$, $C_j$ is
constant for each $j$, and \eqref{eqn14} holds. \end{proof}

\begin{theorem}\label{thmX}
Let $\gamma$ be a curve in ${\cQ}_0$ taking values in $\R^d$. Then
$\gamma$ canonically determines the following data:
\begin{itemize}
\item an orthogonal decomposition of the target space ${\R}^d =
  {\R}^{2n} \oplus {\R}^p$,
\item an invertible skew-symmetric linear map $A$ on $\R^{2n}$,
\item a complex structure $\R^{2n} \simeq \C^n$,
\item vectors $v\in\R^{2n}$ and $w\in\R^p$, for which
\begin{equation}\label{eqn19}
\gamma(s) = \exp(As) v \oplus w.
\end{equation}
\end{itemize}
The curve $\gamma$ is constant if and only if $n=0$ and $d=p$.
It lies in no hyperplane if and only if $d=2n$ and $p=0$.
\end{theorem}

\begin{corollary}\label{cor3}
Let $\mu \in {\cQ}_1$. With $A,v,w$ as above, there are $v_0$ and
$w_0$ such that
\begin{equation}\label{eqn20}
\mu(s) = \left( ( \exp(As) - I) A^{-1} v + v_0 \right) \oplus
\left( ws + w_0 \right).
\end{equation}
\end{corollary}

\begin{proof}
Put $\gamma = D\mu$. Then \eqref{eqn20} follows by integrating
\eqref{eqn19}.
\end{proof}

In Theorem \ref{thmX} and Corollary \ref{cor3} we allow the
possibilities that $n=0$ or $p=0$. Note that $n=0$ if and only if the
curve $\gamma\in{\cQ}_0$ is constant, or equivalently if its
integral $\mu\in{\cQ}_1$ is affine. The copy of ${\R}^{2n}$ is the
{\it horizontal subspace} defined by $\gamma$. The constant vector $w
\in {\R}^p$ is the direction of the vertical component of $\mu$.

To prove Theorem \ref{thmX} we will show that $\gamma$ solves a
constant coefficient ODE. Hence there are constants $\lambda_j \in
{\C}$ and (not identically zero) vector-valued complex polynomials
$q_j$ such that
\begin{equation}\label{eqn18}
\gamma(s) = \sum_{j=1}^K q_j(s) e^{\lambda_j s}.
\end{equation}
We will show that each $\lambda_j$ is purely imaginary and then that
each polynomial $q_j$ is constant. These statements have the following
consequence. By \eqref{eqn14} there is a polynomial $p$ such that
$p(D) \gamma = 0$. The numbers $\lambda_j$ in \eqref{eqn18} are the
roots of $p$. These roots are distinct if and only if each $q_j$ is
constant. 

Conversely, as previously observed, any curve of the form
\eqref{eqn19} is in ${\cQ}_0$ and any curve of the form
\eqref{eqn20} is in ${\cQ}_1$. As a consequence, we characterize step
two Carnot groups in terms of $\cQ_0$ and $\cQ_1$ curves. See Theorem
\ref{thmB} and Corollary \ref{corC}.

\begin{proof}[Proof of Theorem \ref{thmX}]
By the preceding discussion, we
may assume that $\gamma$ is nonconstant. Formula \eqref{eqn18} holds
because $\gamma \in {\cQ}_0$, and $K\ge 1$ because $\gamma$ is
nonconstant. We first claim that all the $\lambda_j$ in \eqref{eqn18}
must be purely imaginary. Write $\lambda_j = \xi_j + i \eta_j$ and
assume that $|\xi_j|\le|\xi_1|$ for all $j$. Then $c_0 =
||\gamma(s)||^2 = e^{2|\xi_1|s} \Phi(s)$ with
$\limsup_{s\to\infty}\Phi(s)>0$. Thus $0=\xi_1$ and hence $0=\xi_j$
for all $j$.

Next we claim that each polynomial $q_j$ is a constant. In other
words, there are no repeated roots (no resonances). The proof is a
similar asymptotic argument. Suppose that the maximal degree of the
$q_j$ is $d$. Then $c_0 = |s|^{2d} \Psi(s)$ with
$\limsup_{s\to\infty}\Psi(s)>0$. Thus $d=0$ and each polynomial $q_j$
is a (nonzero) constant.

To this point we have shown that there is a polynomial $p$ with
distinct, purely imaginary zeros such that $p(D)\gamma = 0$. We write
$q_j = x_j + i y_j$ for vectors $x_j$ and $y_j$ and put $\lambda_j = i
\eta_j$. Since $\gamma$ is real, \eqref{eqn18} yields
\begin{equation}\label{eqn21}
\gamma(s) = \sum_{j=1}^K x_j \cos(\eta_j s) - y_j \sin (\eta_j s).
\end{equation}
The coefficients $x_j$ and $y_j$ in \eqref{eqn21} are elements of ${\R}^d$.

Recall that the $\eta_j$ are distinct. Let $\eta_1,\ldots,\eta_n$
denote the nonzero values and let $\eta_{n+1}=0$ if necessary. Put $v=
\sum_{j=1}^n x_j$ and put $w = \sum_{j > n} x_j$ in \eqref{eqn21}.
Since $\gamma$ is orthogonal to $D\gamma$, we can
rewrite \eqref{eqn21} as
\begin{equation}\label{eqn22}
\gamma(s) = \left( \sum_{j=1}^n x_j \cos(\eta_j s) - y_j \sin(\eta_j
s)\right) \oplus w
\end{equation}
where $\gamma(0) = v \oplus w$.

We claim that the set of vectors $x_1,y_1,\ldots,x_n,y_n$ are linearly
independent. Given the claim their span is a canonical copy of
${\R}^{2n}$ in ${\R}^d$. Thus, given $\gamma \in {\cQ}_0$, we obtain
via \eqref{eqn22} a canonical subspace ${\R}^{2n}$ of ${\R}^d$ and its
orthogonal complement ${\R}^p$.

The linear independence and the existence of the skew-symmetric $A$
follow from the following reasoning. Let $J$ be the usual complex
structure matrix in \eqref{J} and let $A$ denote the direct sum
of the blocks $\eta_j J$ for $j\le n$. Then (for $n\ne 0$) the matrix
$A$ is invertible and $\exp(As)$ is the direct sum of blocks
\begin{equation}\label{eqn24}
\begin{pmatrix} \cos(\eta_j s) & -\sin(\eta_j s) \\ \sin(\eta_j s)  &
  \cos(\eta_j s) \end{pmatrix}.
\end{equation}
By combining \eqref{eqn22} and \eqref{eqn24} one obtains \eqref{eqn19}
for vectors $v$ and $w$. If $n>0$ then $v\ne 0$.

The linear transformation $A$ has $2n$ distinct nonzero eigenvalues
$\pm i\eta_1,\ldots,\pm i\eta_n$. Thus for any nonzero vector $v$,
the vectors $A^jv$ for $0\le j\le 2n-1$ are linearly independent. We
conclude that $\gamma$ lies in no hyperplane if and only if $p=0$.
\end{proof}

Theorem \ref{thmB} and Corollary \ref{corC} provide the decisive
relationship among step two Carnot
groups and their geodesics, helical CR structures, and curves in
$\cQ_0$ and $\cQ_1$. Part (b) of Corollary \ref{corC} restates Theorem
\ref{maintheorem} from the introduction.

\begin{theorem}\label{thmB}
(a) Each curve $\gamma\in\cQ_0$ determines a helical CR structure on
its target space whose horizontal space coincides with that of
$\gamma$. Conversely, each helical CR structure with horizontal
space $H$ determines a family of curves $\gamma\in\cQ_0$,
parameterized by $\pi_H(\gamma(0))$, where $\pi_H$ denotes the
projection from ${\R}^d$ onto $H$.

(b) Each curve $\mu\in\cQ_1$ determines a marked helical CR
structure on its target space whose horizontal space coincides with
that of $\mu$. Conversely, each marked helical CR structure
determines a curve in $\cQ_1$.
\end{theorem}

\begin{corollary}\label{corC}
(a) Each nonconstant curve in $\cQ_0$ contained in a hyperplane of
$\R^d$ determines a unique step two stratified Lie algebra of
contact type on a subspace of $\R^d$. Conversely, each such Lie
algebra with horizontal space $H$ determines a family of curves
$\gamma\in\cQ_0$, parameterized by $\pi_H(\gamma(0))$.

(b) Each nonaffine curve $\mu \in \cQ_1$ contained in a hyperplane
of $\R^d$ determines a unique step two stratified Lie algebra of
contact type on a subspace of $\R^d$ together with the germ of a
normal geodesic $c$ for the induced Carnot-Carath\'eodory metric on
the associated Lie group. Conversely, each such Lie algebra and
geodesic determine a curve in $\cQ_1$. The horizontal projections of
$\mu$ and $c$ coincide.
\end{corollary}

Theorem \ref{thmX} has additional geometric consequences. 
Let $\gamma$
be a curve in $\cQ_0$ which for convenience we assume to be
vertically trivial. By Theorem \ref{thmX} we have $\gamma(s) =
\exp(As)v$ for some skew symmetric $A$ and $v\in\R^{2n}$. Since the
eigenvalues of $A$ are distinct, ${\R}^{2n}$ is the direct sum of
two-dimensional subspaces which by \eqref{J} and \eqref{eqn24} we
may identify with ${\C}$. The projection of $\gamma$ into each of
these particular planes is a circle. Of course $\gamma$ is in
general not a circle, because its image does not lie in a plane. A
well-known example is the {\it skew-line} $\phi:\R\to
\Sph^1\times\Sph^1$ defined by $\phi(s) = (e^{is}, e^{i \alpha s})$
for irrational $\alpha$, which embeds $\R$ injectively in
$\Sph^1\times\Sph^1$. Let $\gamma$ be an element of ${\cQ}_0$ which
maps into no hyperplane in ${\R}^{2n}$. Theorem \ref{thmX} allows us
to identify ${\R}^{2n}$ with ${\C}^n$ and to think of $\gamma$ as
being defined by
\begin{equation}\label{eqn25}
\gamma(s) = (\zeta_1 e^{i \alpha_1 s},\ldots,\zeta_n e^{i \alpha_n s}).
\end{equation}
From \eqref{eqn25} we obtain a simple criterion for injectivity of
$\cQ_1$ curves. 

\begin{corollary}\label{cor4}
A curve $\gamma\in\cQ_1$ fails to be injective if and only if each
$\alpha_j$ in \eqref{eqn25} is a rational multiple of $\alpha_1$. In
this case $\gamma$ traces its image countably many times.
\end{corollary}

\begin{proof}
Injectivity fails if and only if there are distinct
real $s$ and $u$ such that $e^{i \alpha_j u} = e^{i \alpha_j s}$ for
all $j$, equivalently, if and only if $\alpha_j (u-s) = 2\pi k_j$
for some integers $k_j$.
\end{proof}

\subsection{Homogeneous elements of ${\cQ}_0$}\label{sec4}
In this section we study specific curves in $\cQ_0$ defined by
homogeneous polynomial expressions analogous to those in \cite{D1}.
The study of these examples motivated the present work. We describe
the associated helical CR 
structures and determine the spectra of the relevant skew-symmetric
matrices. In Remark \ref{motivate} we briefly discuss the
corresponding polynomial mappings in several complex variables.

For $0 \le j \le m$ let $E_j$ denote the standard $j$th basis vector
of ${\R}^{m+1}$. We define the curve $\gamma_m:{\R} \to {\R}^{m+1}$ by
the formula
\begin{equation}\label{eqn27}
\gamma_m(s) = \sum_{j=0}^m {\sqrt{m \choose j}} \cos^{m-j}(s) \sin^j(s)
E_j.
\end{equation}
We say that a curve $\gamma$ with values in ${\R}^d$ is {\it
homogeneous of degree} $m$ if
\begin{equation}\label{eqn26}
\gamma(s) = H(\cos s,\sin s)
\end{equation}
for some $H \in V_m(2,d)$, where $V_m(2,d)$ denotes the collection of
homogeneous polynomial maps of degree $m$ in two real variables with
values in ${\R}^d$.

\begin{proposition}\label{prop3}
For each $m$, $\gamma_m$ lies in the unit sphere,
$\gamma_m \in {\cQ}_0$, and $\gamma_m$ is homogeneous of degree $m$.
\end{proposition}

\begin{proof}
From \eqref{eqn27} we see that $\gamma_m$ is homogeneous of degree
$m$ and that $\gamma_m(0)$ is the unit vector $E_0$. We claim that
$\gamma_m$ satisfies the ODE $D\gamma_m = L_m \gamma_m$ where $L_m$ is
skew-symmetric. It follows that
$||\gamma_m(s)|| =  ||\gamma_m(0)|| =1$ for all $s$. Thus it suffices
to prove the claim.

We determine $L_m$ by using the elementary identity
\begin{equation}\label{eqn28}
\frac{d}{ds} \left(\cos^a(s)  \sin^b (s) \right)
 = b \cos^{a+1}(s)  \sin^{b-1}(s)  - a \cos^{a-1}(s) \sin^{b+1}(s).
\end{equation}
With respect to the standard basis the matrix entries of $L_m$
all vanish except on the super- and sub-diagonals. Furthermore each
entry on the superdiagonal is negative, and the corresponding entry
on the subdiagonal is its additive inverse. It follows that $L_m$ is
skew-symmetric and hence that $\gamma_m \in {\cQ}_0$.
\end{proof}

We give a second explanation for why $\gamma_m \in {\cQ}_0$. By part
(iii) of Lemma \ref{lem1}, the tensor product of elements in ${\cQ}_0$
is also in ${\cQ}_0$. We can identify $\gamma_m$ with the $m$-fold
symmetric tensor product of $\gamma_1$, whose image is the unit circle.
It is also easy to see why $\gamma_m$ lies in the unit sphere.
By \eqref{eqn9} we have
\begin{equation}\label{eqn29}
||\gamma_m (s)||^2 =  ||\gamma_1 (s)||^{2m} = \left( \cos^2 s +
\sin^2 s \right)^m = 1.
\end{equation}
Let $\xi$ be a curve that is homogeneous of degree $m$. The monomials
$x^{m-k} y^k$ for $0 \le k \le m$ form a basis for $V_m(2,1)$.
The component functions of $\gamma_m$ therefore span the space of
homogeneous polynomials of degree $m$ in $\cos s$ and $\sin s$. It
follows that there a linear mapping $B$ such that
\begin{equation}\label{eqn30}
\xi(s) = B \gamma_m (s) = B \exp(L_m s) E_0.
\end{equation}
The mapping $B$ need not be orthogonal, nor even invertible, even if
$\xi \in {\cQ}_0$.

\begin{example}\label{ex2}
Let $ \xi(s) = (\cos^2 s - \sin^2 s, 2 \cos s \sin s, 0)$. Then $\xi$
is homogeneous of degree $2$, and hence can be written $B\gamma_2$ for
some $B$. In this case $B$ is not invertible, yet $\xi \in
{\cQ}_0$. The image is the unit circle.
\end{example}

We write matrix representations for the first few $L_m$ and explicit
formulas for the corresponding $\gamma_m$. By convention $\gamma_0$
is the constant $1$. We have
$$
L_1 = \begin{pmatrix} 0 & -1 \\  1 & 0\end{pmatrix},
$$
\begin{equation}\label{eqn31}
\gamma_1(s) = (\cos s, \sin s),
\end{equation}
$$
L_2 = \begin{pmatrix} 0 & - \sqrt{2} & 0 \\  \sqrt{2} & 0 & - \sqrt{2}
  \\ 0 & \sqrt{2} & 0 \end{pmatrix},
$$
\begin{equation}\label{eqn32}
\gamma_2(s) = (\cos^2 s, \sqrt{2} \cos s \sin s, \sin^2 s ),
\end{equation}
$$
L_3 = \begin{pmatrix} 0 & - \sqrt{3} & 0 & 0 \\ \sqrt{3} & 0 & - 2 & 0
  \cr 0 & 2 & 0 & - \sqrt{3} \\ 0 & 0 & \sqrt{3} & 0 \end{pmatrix},
$$
\begin{equation}\label{eqn33}
\gamma_3(s) = (\cos^3 s, \sqrt{3} \cos^2 s\sin s, \sqrt{3} \cos
s\sin^2 s, \sin^3 s ).
\end{equation}
Since $L_m$ is skew-symmetric, all its eigenvalues are purely
imaginary. In fact, for each $m$, the eigenvalues of $L_m$ are
distinct, and hence $L_m$ is diagonalizable over ${\C}$. We next
determine the characteristic polynomial of the $(m+1)\times(m+1)$
matrix $L_m$. We write $p_m(x)$ for the characteristic polynomial
$\det(L_m - x I)$.

\begin{proposition}\label{prop4}
The characteristic polynomial $p_m(x)$ is
$$ p_{2k} (x) = - x \prod_{j=1}^k (x^2 + (2j)^2) $$
if $m=2k$ is even, and
$$ p_{2k+1} (x) = \prod_{j=0}^k (x^2 + (2j+1)^2) $$
if $m=2k+1$ is odd. The eigenvalues of $L_{2k}$ are $ 0, \pm 2i, \pm
4i, \ldots, \pm ki$, while the eigenvalues of $L_{2k+1}$ are $ \pm i,
\pm 3i, \pm 5i,\ldots, \pm (2k+1) i$.
\end{proposition}

We denote by $\sigma(A)$ the spectrum of an operator $A$, and by
$S+T$, resp.\ $S\cdot T$, the Minkowski sum $\{s+t:s\in S,t\in T\}$,
resp.\ Minkowski product $\{s\cdot t:s\in S,t\in T\}$, of two sets
$S,T\subset\C$. We write $m\cdot S = \underbrace{S+\cdots+S}_m$ and
$S^m = \underbrace{S\cdot \cdots \cdot S}_m$.

\begin{proof}
Recall that $\gamma_m$ can be identified with the $m$-fold (symmetric)
tensor product of $\gamma_1$. We have $\gamma_1(s) = {\rm exp}( L_1
s)e_0 = {\rm exp}(Js) e_0,$ where $e_0=(1,0)$ in ${\R}^2$. Therefore
$$
\gamma_m(s) = \left( {\exp}(Js)e_0 \right)^{\otimes m} = ({\exp}(Js))^{\otimes
m}E_0 = {\exp}(L_m s)E_0.
$$
For operators $A,B$, $\sigma(A \otimes B) = \sigma(A)\cdot\sigma(B)$.
Thus
$$
\sigma(\exp(L_m s)) = \sigma(\exp(Js))^m = \{e^{is},e^{-is}\}^m
$$
and $\sigma(L_m) = m \cdot \{ \pm i \} = \{ \pm
(m-2j)i : j=0,1,\ldots,[\frac{m}{2}]\}$, where $[\cdot]$ denotes the
greatest integer function. The formulas for the
characteristic polynomials follow.
\end{proof}

We mention a related suggestive point of view. Let $I_{m+1}$ denote
the identity operator on ${\C}^{m+1}$. Define a linear operator $
{\cL}_m$ by
\begin{equation}\label{eqn34}
{\cL}_m = \frac1\pi\log( (-1)^m I_{m+1}).
\end{equation}
This suggestive notation means that $ {\cL}_m $ is the diagonal matrix
whose eigenvalues are $\frac1\pi$ times $m+1$ particular values of
$\log((-1)^m)$. When $m=2k$ these values are $$ mi,
(m-2)i,\ldots,2i,0, -2i,\ldots, -(m-2)i, -m i  $$ and when $m=2k+1$
they are $$ mi, (m-2)i,\ldots,i, -i,\ldots, -(m-2)i, -m i. $$ The
operators ${\cL}_m$ and $L_m$ are similar over ${\C}$ and hence have
the same eigenvalues.

We observe explicitly part of Theorem \ref{thmX}. When $m+1$ is odd,
one of the eigenvalues is zero. It follows in this case that
$\gamma_m$ maps to a hyperplane. For example, the image of $\gamma_2$,
a priori in ${\R}^3$, is in fact a circle in the hyperplane
$\{(x_0,x_1,x_2):x_0+x_2=1\}$, similarly, the image of $\gamma_4$ lies
in a four-dimensional hyperplane in ${\R}^5$.

\begin{remark}\label{motivate}
We mention another role played by the curves $\gamma_m$. Using
homogeneity and polar coordinates we extend $\gamma_m$ to all of
${\R}^2$ by setting
$$ P_m(r\cos\theta,r\sin\theta) = r^m \gamma_m(\cos \theta, \sin
\theta).
$$
Then $P_m$ defines a proper polynomial mapping from the unit disk to
the unit ball in ${\R}^{m+1}$. This map $P_m$ is invariant under a
cyclic subgroup of order $m$ of $O(2)$. See \cite{D4} and its
references for information about holomorphic mappings invariant
under finite subgroups of the unitary group. The particular
group-invariant mappings $z \mapsto z^{\otimes m}$ are the simplest
examples, and their restrictions to the unit sphere in $\C^m$
provide a complex variables analogue of the curves $\gamma_m$.
\end{remark}

\section{Concluding Remarks}

\begin{remark}\label{rem1}
We discuss the issue of reparameterization. In this paper we have
considered curves as maps from an interval in ${\R}$ into ${\R}^d$. As
previously mentioned, the class of curves in ${\cQ}_0$ is closed under
affine reparameterization $s\mapsto \lambda s + b$. Under the
correspondence in Theorem \ref{thmB}, such reparameterizations
correspond to an equivalence of helical CR structures, where two such
structures ${\cC}$ and ${\cC}'$ are called {\it equivalent} if the
associated skew-symmetric matrices $A$ and $A'$ satisfy $A'=\lambda A$
for some nonzero $\lambda$. Under the correspondence in Corollary
\ref{corC}, such reparameterizations correspond to isomorphic
stratified Lie algebra structures. Similar statements can be made for
curves in ${\cQ}_1$; we leave the precise statements to the reader.
\end{remark}

\begin{remark}\label{DRZ}
The skew-symmetric matrices associated to the curves $\gamma_m$ 
are bidiagonal, i.e., have
non-zero entries only on the super- and sub-diagonals. Closely related
issues arise in the fundamental algebraic question of normal forms for
orthogonal similarity classes of skew-symmetric matrices. The paper
\cite{DRZ} considers a skew-symmetric matrix $A$ defined over an
algebraically closed field ${\F}$ and seeks an orthogonal $P$ such that
$PAP^{-1}$ is as close to bidiagonal as possible. Achieving
bidiagonality itself is not always possible. In our context the
underlying skew-symmetric matrices are bidiagonal and real.  
\end{remark}

\begin{remark}\label{rem2}
We can extend Corollary \ref{corC} to relate collections of curves
in ${\cQ}_1$ to sub-Riemannian structures with additional vertical
directions, enabling us to characterize arbitrary step two Carnot
groups. For each 
$\alpha=1,\ldots,p$, let $\mu_\alpha$  be a nonaffine curve in
${\cQ}_1$ which lies in a hyperplane in ${\R}^d$; we assume that the
horizontal spaces of these curves coincide. We further assume that
the vertical directions $w_\alpha$ for the curves $\mu_\alpha$ are
linearly independent in the orthogonal complement ${\R}^p$. A
straightforward induction using Corollary \ref{corC} leads from this
data to a unique step two stratified Lie algebra of type $(2n,p)$
together with $p$ normal geodesics $c_1,\ldots,c_p$ for the CC
metric on the associated Lie group. Conversely, beginning from such
a Lie algebra together with $p$ geodesics, we construct $p$ curves
$\mu_1,\ldots,\mu_p$ in ${\cQ}_1$ taking values in ${\R}^d$ with
common horizontal space and vertical directions $w_1,\ldots,w_p$
forming a basis for ${\R}^p$. We obtain the following theorem.

\begin{theorem}\label{thm4}
Each $p$-tuple of nonaffine curves in ${\cQ}_1$, each of which lies
in a hyperplane in ${\R}^d$, with common horizontal space, and with
vertical directions $w_1,\ldots,w_p$ forming a basis for ${\R}^p$,
uniquely determines and is uniquely determined by a step two
stratified Lie algebra ${\bf g}$ on ${\R}^{d}$ together with a
$p$-tuple of distinguished normal geodesics for the CC metric on the
associated Carnot group.
\end{theorem}
\end{remark}


\begin{thebibliography}{MM}
\bibitem{B1} Roger W. Brockett, {\it Nonlinear control theory and
    differential geometry} in ``Proceedings of the International
    Congress of Mathematicians'', Vol.\ 1,2 (Warsaw, 1983), 1984,
    1357--1368.
\bibitem{B2} Roger W. Brockett, {\it Control theory and singular Riemannian
geometry} in ``New directions in applied mathematics'', 1982, 11--27.
\bibitem{B} Andr{\'e} Bella{\"{\i}}che, {\it The tangent space in
sub-Riemannian geometry} in ``Sub-Riemannian Geometry'', Progress in
Math., Birkh{\"a}user, 144, 1996, 1--78.
\bibitem{BG} Thomas Bloom and Ian Graham, {\it On `type' conditions
    for generic real submanifolds of $\C^n$}, Invent.\ Math.\ 40
  (1977), no.\ 3, 217--243.
\bibitem{CDPT} Luca Capogna, Donatella Danielli, Scott Pauls, and Jeremy T.
Tyson, ``An Introduction to the Heisenberg group and the
sub-Riemannian isoperimetric problem'', Progress in Math.,
Birkh{\"a}user, 259, 2007.
\bibitem{CL} Luca Capogna and Fang-Hua Lin, {\it Legendrian energy
    minimizers. {I}. {H}eisenberg group target}, {Calc.\ Var.\ PDE},
  12 (2001), 145--171.
\bibitem{D1} John P. D'Angelo, {\it A monotonicity result for volumes of
holomorphic images}, Michigan Math J., 54 (2006), 1--24.
\bibitem{D2} John P. D'Angelo, ``Several Complex Variables and the
  Geometry of Real Hypersurfaces'', CRC Press, Boca Raton, FL, 1992.
\bibitem{D4} John P. D'Angelo, {\it Invariant holomorphic maps}, J.
  Geom.\ Anal.\ 6 (1996), 163--179.
\bibitem{DLP} John P. D'Angelo, Jiri Lebl, and Han Peters, {\it Degree
estimates for polynomials constant on a hyperplane}, Michigan Math.
J., 55 (2007), no.\ 3, 693--713.
\bibitem{DRZ} Dragomir Z. Dokovic, Konstanze Rietsch and Kaiming Zhao, {\it
Normal forms for orthogonality similarity classes of skew-symmetric
matrices}, preprint, 2006, arXiv:math.RT/0603245 v2.
\bibitem{Ga} Bernard Gaveau, {\it Principe de moindre action, propagation de
la chaleur et estim\'ees sous elliptiques sur certains groupes
nilpotents}, Acta Math. 139 (1977), 95--153.
\bibitem{Gl1} Herman Gluck, {\it Higher curvatures of curves in
    Euclidean space}, American Math Monthly 73 (1966), 699--704.
\bibitem{Gl2} Herman Gluck, {\it Higher curvatures of curves in
    Euclidean space II}, American Math Monthly 74 (1967), 1049--1056.
\bibitem{gk} Ch.\ Gol\'e and R. Karidi, {\it A note on {C}arnot
    geodesics in nilpotent {L}ie groups}, J. Dynam.\ Control Systems,
  1 (1995), no.\ 4, 535--549.
\bibitem{M} Richard Montgomery, ``A tour of subriemannian geometries, their
geodesics and applications'', Mathematical Surveys and Monographs 91,
A.M.S. 2002.
\bibitem{RS} Linda Preiss Rothschild and E. M. Stein,
{\it Hypoelliptic differential operators and nilpotent groups},
Acta Math., 137 (1976), no.\ 3-4, 247--320.
\bibitem{S1} Robert S. Strichartz, {\it Sub-Riemannian geometry}, J.
  Differential Geom., 24 (1986), 221--263.
\bibitem{S2} Robert S. Strichartz, {\it Corrections to
    ``Sub-Riemannian geometry'' [J. Differential Geom., 24 (1986),
    221--263]}, J. Differential Geom., 30 (1989), 595--596.
\end{thebibliography}
\end{document}